\documentclass[12pt]{article}
\usepackage{mathrsfs}
\usepackage{amssymb}
\usepackage[tbtags]{amsmath}
\usepackage{cite}
\usepackage[noblocks]{authblk}
\numberwithin{equation}{section}

\newcommand{\ba}{\begin{array}}
\newcommand{\ea}{\end{array}}

\newcommand{\xiaowuhao}{\fontsize{8pt}{\baselineskip}\selectfont}

\newcommand {\xo}{\Omega}

\date{}
\begin{document}
\title{\textbf{Existence and regularity results for anisotropic parabolic equations with
degenerate coercivity}}
 \author{Weilin Zou$^{\mbox{a,}}$\footnote{\ Corresponding author.\\
\ \ \ \ \indent \ \ \ \ \ E-mail addresses: zwl267@163.com, Tel.: +8679183863755.
},~~~~Yuanchun Ren$^{\mbox{a}}$, ~~~~Wei Wang$^{\mbox{b}}$ \\
\scriptsize{ \xiaowuhao $^{\mbox{a}}$ College of Mathematics and Information
Science, Nanchang Hangkong University, Nanchang, 330063, China}\\
\scriptsize{ \xiaowuhao $^{\mbox{b}}$ Party and Administration Office, Nanchang Hangkong University, Nanchang, 330063, China}} \maketitle

 {\bf Abstract:} This paper deals with a class of nonlinear anisotropic parabolic equations with  degenerate coercivity. Using the anisotropic Gagliardo-Nirenberg-type inequality, we prove some existence and regularity results for the solutions under the framework of anisotropic Sobolev spaces, which generalize the previous results of \cite{f,d,e}.



\section{Introduction and statement of the main results}
\quad This paper deals with the existence and regularity results of weak solutions to the following anisotropic parabolic problem
 $$
(P)\ \left\{
\begin{array}{ll}
\frac{\partial u}{\partial t}-\overset{N}{\underset{i=1}\sum} D_i(a_i(x,t,u,Du)) =f & \mbox{in } Q,\\ u=0 &
\mbox{on }\partial \Omega\times (0, T),\\
 u(x, 0)= 0 & \mbox{in }\Omega,
\end{array}
\right. $$
where $f\in L^m(Q)$ with $m\geq 1$,
$\Omega$ is an open bounded subset of $\mathbb{R}^N(N\ge 2)$ , $T$ is a positive constant
and $Q=\Omega\times (0, T)$ with the lateral boundary
$\partial\Omega\times (0, T)$.

Here $a_i: Q\times \mathbb{R}\times \mathbb{R}^N \rightarrow \mathbb{R}$ ($i=1,2...,N$) are Carath\'{e}odory
functions such that for almost all $(x,t)\in Q$, for all $s\in
\mathbb{R}$,  and for all $\xi,\zeta\in\mathbb{R}^N$ with $\xi\neq\zeta$,
\begin{align}\label{1.2}
\sum_{i=1}^{N}a_i(x,t,s,\xi)\xi_i\geq \alpha\sum_{i=1}^{N}\frac{|\xi_i|^{p_i}}{(1+|s|)^{\theta(p_i-1)}},
\end{align}
\begin{align}\label{1.1}0\leq\theta<1+\frac{\bar{p}}{N(\bar{p}-1)}, \end{align}
\begin{align}\label{1.3}
|a_i(x,t,s,\xi)|\leq\beta(1+|\xi_i|^{p_i-1}),\ \ \ i=1,2,\ldots ,N,
\end{align}
and
\begin{align}\label{1.4}
\sum_{i=1}^{N} [ a_i(x,t,s,\xi )-a_i(x,t,s,\zeta ) ][ \xi_i -\zeta_i ]>0,
\end{align}
where  $\alpha, \beta$
are two positive constants, $p_1,...,p_N \in (1,+\infty)$ and $\bar{p}$ is defined as in (\ref{2.001}).

As $\theta=0$, the typical model of problem $(P)$ reduces to the so-called anisotropic evolutionary $\vec{p}$-Laplacian equation. Indeed, assumption (\ref{1.2}) presents an interesting feature, namely an anisotropic diffusion with orthotropic structure. These kinds
of problems have attracted an increasing amount of attention in recent years as they reflect the modeling of physical and mechanical processes in anisotropic continuous medium. With no hope of being thorough, we mention
some papers regarding the study of these problems \cite{ming,Feo21,Barletta2020,G,42,43,32,33,34,f,h} (see also the references therein).

In the case that $\theta\neq0$, the parabolic operator in problem $(P)$ may degenerate as soon as the solution is unbounded. Indeed, the
diffusion coefficient may tend to zero as the solution $u$ becomes large, which means that a slow diffusion effect may appear. There
are already some results for such kind of problems both in the
stationary and evolution cases with isotropic growth condition. For example, the isotropic elliptic case (with $p_i\equiv p$ for $i=1,2,..., N$) associated to problem $(P)$ was firstly studied in \cite{21} and then developed in \cite{bo03,B, D, E, Zou21}; in the isotropic parabolic case, the existence and regularity results of problem $(P)$ were proved in \cite{d}(see also \cite{a,b,por14}), and the study of problem $(P)$ with nonzero initial data can be found in \cite{por08}.

Recently, the regularity results for weak solutions to the anisotropic elliptic case of problem $(P)$ were investigated in \cite{e}. The existence and regularity of the entropy solutions for stationary problem $(P)$ with lower order terms were considered in \cite{ayadi20}. However, to the best of our knowledge, there are no existence and regularity results concerning with the nonlinear anisotropic parabolic equations with degenerate coercivity. Hence, the aim of this paper is to extend the results of \cite{d,e} to the case of anisotropic parabolic problems. To deal with this kind of problem, we shall look for $u_n$ as a solution to the approximate problem with non-degenerate coercivity, and establish some estimates on the weak solutions of these problems under the framework of anisotropic Sobolev spaces. The proofs of these estimates rely crucially on the anisotropic version of Gagliardo-Nirenberg inequality and the interpolation inequality.

The remainder of this paper is
organized as follows.  In Section 2, we give some preliminary results and state the main results. Section 3 is devoted to proving some estimates of weak solutions to the approximate problem under different assumptions on $f$. In Section 4, we briefly give the proofs of Theorem 2.1-2.5.

\section{Preliminaries and main results}

$\qquad$ Given $\vec{p}=(p_1,p_2,...,p_N)$, we define
 \begin{equation}\label{2.001}
 \bar{p}=\left(\frac{1}{N}\sum\limits_{i=1}^{N}\frac{1}{p_i}\right)^{-1}.\end{equation}
Assume $1<\bar{p}<N$, we denote
 $$\bar{p}^*=\frac{N\bar{p}}{N-\bar{p}}.$$

{\bf Anisotropic Sobolev space.} Let $\Omega$ be an open bounded subset of $\mathbb{R}^N(N\ge 2)$, we define the following space
$$ W_0^{1,\vec{p}}(\Omega):=\left\{g\in W_0^{1,1}(\Omega):D_i g\in L^{p_i}(\Omega)\ \mbox{for\ every}\ i=1,2\ldots N\right\},$$
which is a Banach space under the norm
$$ \|g\|_{W_0^{1,\vec{p}}(\Omega)}:=\sum\limits_{i=1}^{N}\|D_ig\|_{L^{p_i}(\Omega)}.$$
We also introduce the functional space $L^{\vec{p}}(0,T;W_0^{1,\vec{p}}(\Omega))$ by
$$L^{\vec{p}}(0,T;W_0^{1,\vec{p}}(\Omega)):=\left\{u\in L^{1}\left(0,T;W_0^{1,1}(\Omega)\right):D_i u\in L^{p_i}(Q),\ \mbox{for\ every}\ i=1,2,\ldots ,N\right\}.$$
The dual space of $L^{\vec{p}}(0,T;W_0^{1,\vec{p}}(\Omega))$ is denoted by $L^{\vec{p'}}(0,T;W^{-1,\vec{p'}}(\Omega))$, where $\vec{p'}=(p_1',p_2',...,p_N')$.

Next we recall the following anisotropic Sobolev embedding Theorem(see \cite{26,27}).\vspace{2mm}

\noindent{\bf Lemma 2.1.} \emph{Let $u\in W_0^{1,\vec{p}}(\Omega)$. If $1<\bar{p}<N$, then there exists a constant $C = C(N, \vec{p}) > 0$ such that
$$\|u\|_{L^{\bar{p}^*}(\Omega)}\leq C\left(\prod_{i=1}^{N}\|D_iu\|_{L^{p_i}(\Omega)}\right)^{1/N}.$$}\vspace{2mm}

Now we give the definition of weak solutions of problem $(P)$.\vspace{2mm}

\noindent{\bf Definition 2.1.} A measurable function $u \in
L^1(0,T; W_0^{1,1}(\Omega))$ is a weak solution of
problem $(P)$ if $a_i(x,t,u,Du)\in L^1(Q)$ and
 \begin{equation}
 \label{7}-\iint_Q
u\phi_tdxdt +\iint_Q \sum_{i=1}^{N}a_i(x,t,u,Du)D_i\phi dxdt=\iint_Q f\phi dxdt,\notag
\end{equation}
for every test function $\phi\in C^\infty(\bar{Q})$ such that $\phi=0$ in a neighborhood of $\partial\Omega\times(0,T)\cup (\Omega\times \{T\})$.\vspace{2mm}

The main results of this work are stated as follows.\vspace{2mm}

\noindent{\bf Theorem 2.1.} \emph{ Assume that (\ref{1.2})-(\ref{1.4}) hold and $f\in
L^m(Q)$ with $m>\frac{N}{\bar{p}}+1$. Then there exists
 a weak solution $u\in L^{\vec{p}}(0,T;W_0^{1,\vec{p}}(\Omega))\cap
L^\infty(Q)$ to problem $(P)$.}\vspace{2mm}

\noindent{\bf Remark 2.1.} This result is independent of $\theta$ which corresponds to the one that obtains in the isotropic case with $\theta=0$(see \cite{7}). If $p_1=p_2=...=p_N=2$, it is proved in \cite{d} that problem $(P)$ admits at least a bounded weak solution provided that $f$ belongs to $
L^m(Q)$ with $m>\frac{N}{2}+1$. Thus, Theorem 2.1 extends the results of \cite{7} and \cite{d}.\vspace{2mm}

\noindent{\bf Theorem 2.2.}  \emph{Assume that (\ref{1.2})-(\ref{1.4}) hold and $f\in
L^m(Q)$ with $m=\frac{N}{\bar{p}}+1$. Then there
exists a constant $\lambda>0$ such that problem $(P)$ admits at least a weak solution
$u\in L^{\vec{p}}(0,T;W_0^{1,\vec{p}}(\Omega))$ with $e^{\lambda |u|^{1-\frac{N\theta(\bar{p}-1)}{N(\bar{p}-1)+\bar{p}}}}\in L^1(Q)$.}\vspace{2mm}

\noindent{\bf Remark 2.2.} It is clear that $e^{\lambda |u|^{1-\frac{N\theta(\bar{p}-1)}{N(\bar{p}-1)+\bar{p}}}}\in L^1(Q)$ implies
$u\in L^{\tilde{r}}(Q)$ for any $\tilde{r}\in (1,+\infty)$. This improves the result of Theorem 1.2 in \cite{d}.\vspace{2mm}
%$m=\frac{N}{\bar{p}}+1$ is shown by Theorem 1.2. For what I can tell, I
%fail in finding other papers which figure out the limit case for parabolic
%equations even though $\theta=0.$

\noindent{\bf Theorem 2.3.} \emph{Assume that (\ref{1.2})-(\ref{1.4}) and $f\in
L^m(Q)$ with $m$ satisfying
\begin{equation}
\label{1.03}m_1\leq m<\frac{N}{\bar{p}}+1,
\end{equation}
where
\begin{equation}\label{0204}
m_1:=\frac{N+2\bar{p}+N(\bar{p}-1)(1-\theta)+\theta\left[\left(\max\limits_{1\leq i\leq N}p_i-1\right)(N+\bar{p})\right]}{N(1-\theta)(\bar{p}-1)+2\bar{p}-\theta\bar{p}\left(\max\limits_{1\leq i\leq N}p_i-1\right)}.
\end{equation}
Then there exists
 a weak solution $u\in L^{\vec{p}}(0,T;W_0^{1,\vec{p}}(\Omega))\cap
L^{\tilde{r}}(Q)$ problem $(P)$ with
\begin{equation}
 \label{1.04}\tilde{r}=\frac{mN[(1-\theta)(\bar{p}-1)]+m\bar{p}}{N-m\bar{p}+\bar{p}}.
\end{equation}}\vspace{2mm}

\noindent{\bf Theorem 2.4.} \emph{Assume that (\ref{1.2})-(\ref{1.4}) hold and $f\in
L^m(Q)$ with $m$ satisfying
\begin{equation}\begin{split}
 \label{1.5}\frac{N+2+\theta(\bar{p}-1)}{\bar{p}(1+N)+(1-N)[1+\theta(\bar{p}-1)]}<m<m_N
\end{split},\end{equation}
where
\begin{equation}\label{206}
m_N:=\frac{N+2\bar{p}+N(\bar{p}-1)(1-\theta)+\theta\left[\left(\min\limits_{1\leq i\leq N}p_i-1\right)(N+\bar{p})\right]}{N(1-\theta)(\bar{p}-1)+2\bar{p}-\theta\bar{p}\left(\min\limits_{1\leq i\leq N}p_i-1\right)}.
\end{equation}
 Then there exists a weak solution $u\in L^{\vec{q}}(0,T;W_0^{1,\vec{q}}(\Omega))\cap
L^{\tilde{r}}(Q)$ to problem $(P)$ with $\vec{q}=(q_1,q_2,...,q_N)$, where $\tilde{r}$ is defined as in (\ref{1.04}) and
\begin{equation}\label{2.08}
q_i=\frac{m p_i[N(1-\theta)(\bar{p}-1)+\bar{p}]}{[N-\bar{p}(m-1)][\theta(p_i-1)+1]+N(\bar{p}-1)(1-\theta)+\bar{p}}.
\end{equation}}

\noindent{\bf Remark 2.3.} The case of $f\in L^m(\Omega)$ with $m_N\leq m<m_1$ is not considered in Theorem 2.1- Theorem 2.4. This case is rather more complicate, which will be discussed in Remark 4.1 of Section 4. If $1\leq m<\frac{N+2+\theta(\bar{p}-1)}{\bar{p}(1+N)+(1-N)[1+\theta(\bar{p}-1)]}$, it is possible to prove the existence of an entropy solution in the sense made precisely below(see
\cite{12,30,15,16,17,18}).\vspace{2mm}

For this we need to introduce two usual functions defined as follows
\begin{equation}
\begin{split}\label{1.8}
 T_k(s)=\min\{k,\max\{-k,s\}\}, \ \ S_k(s)=\int^s_0 T_k(\tau)d\tau,\ \ \forall s\in \mathbb{R}.
\end{split}
\end{equation}

\noindent{\bf Definition 2.2.} A measurable function  $u \in
L^\infty(0,T; L^1(\Omega))$  is called an entropy solution of
problem $(P)$, if it satisfies the following two conditions:
\begin{equation}\label{210}T_k(u)\in L^{\vec{p}}(0,T;W_0^{1,\vec{p}}(\Omega)),\end{equation}
\begin{equation}\label{211}S_k(u)\in C([0,T];L^1(\Omega)),\end{equation}
and if
\begin{equation}
\begin{split}\label{2.2}&\int_{\Omega}S_k\left(u(T)-\phi(T)\right)dx-\int_{\Omega}S_k\left(-\phi(0)\right)dx
+\int_0^T \left\langle\phi_t,
T_k(u-\phi)\right\rangle dt\\ & +\iint_Q\sum_{i=1}^{N} a_i(x,t,u,Du)D_iT_k(u-\phi) dxdt= \iint_Q fT_k(u-\phi) dxdt,
\end{split}
\end{equation}
 for every
$\phi\in L^{\vec{p}}(0,T;W_0^{1,\vec{p}}(\Omega))\cap L^\infty(Q)$ such that
$\phi_t\in L^{\vec{p'}}(0, T; W^{-1,\vec{p'}}(\Omega))\cap L^{1}(Q).$\vspace{2mm}

\noindent{\bf Remark 2.4}  It is well known that conditions (\ref{210}) and (\ref{211}) allow to define the very weak gradient $Du$ almost everywhere in $Q$(see \cite{12}). For any $k>0$, we have $$ DT_k(u)=v\chi_{\{|u|<k\}}\  \
\mbox{a.e. in }\ \ Q,$$ where $\chi_{\{|u|<k\}}$ represents the
characteristic function on the set $\{|u|<k\}$.\vspace{2mm}

\noindent{\bf Definition 2.3.} (see \cite{30,19})  We say that a measurable function $u:Q\rightarrow R$ belongs to the Marcinkiewicz space ${\cal M}^q(Q)$ with $0<q<+\infty$, if there exists a constant $C$ such that
 $$[u]_q=\mathop {\mbox{sup}}\limits_{k>0}k \ \mbox {meas}\{(x,t)\in Q:
|u(x,t)|>k\}^{\frac{1}{q}}\leq C.$$\vspace{2mm}

For $r<q,$ one can deduce that ${\cal M}^q(Q)\subset {\cal M}^r(Q)$
 and $L^q(Q)\subset {\cal M}^q(Q)\subset L^r(Q)$(see \cite{19,20}).\vspace{2mm}

\noindent{\bf Theorem 2.5.} \emph{Assume that (\ref{1.2})-(\ref{1.4}) hold and $f\in
L^m(Q)$ with
\begin{equation}
\begin{split}\label{1.9}
1\leq m\leq \max \left\{\frac{N+2+\theta(\bar{p}-1)}{\bar{p}(1+N)+(1-N)[1+\theta(\bar{p}-1)]},1\right\}.
\end{split}
\end{equation}
Then there exists an entropy solution $u$ to problem $(P)$ with $u\in\mathcal{M}^{\tilde{r}}(Q)$ and $|D_iu|\in\mathcal{M}^{q_i}(Q)$ for $i=1,2\ldots N$, where $\tilde{r}$ and $q_i$ are defined as before. }

\noindent{\bf Remark 2.5.} If $\theta=0$, Theorem 2.5 reduces to the corresponding result obtained in \cite{11} for parabolic equations with measure data, and
\cite{f} for anisotropic parabolic equations with measure data. \vspace{2mm}

\section{The approximating problems and some useful estimates}

For convenience, we shall use $|E|$ denoting the Lebesgue measure of measurable set $E$. If not otherwise specified, we will denote by $C$ or $C_i$(i=1,2\ldots)several positive constants whose values may change from line to
line. These values will only depend on the data but they will never depend
on the indexes $n$ of the sequences.

To prove Theorem 2.1-2.5, first of all, we consider the following approximate
 problem

 $$
(P_n)\ \left\{
\begin{array}{ll}
\frac{\partial u_n}{\partial t}-\sum\limits_{i=1}^{N}D_i\left(a_i(x,t,T_n(u_n),Du_n\right) =f_n & \mbox{in } Q,\\
u_n=0 & \mbox{on }\partial \Omega\times (0,T),\\ u_n(x,0)=0&
\mbox{in } \Omega,
\end{array}
\right. $$ where $f_n\in C_0^{\infty}(Q)$ such that
 \begin{equation}
\label{2.3}||f_n||_{L^m{(Q)}}\le ||f||_{L^m{(Q)}},\ \ \ \forall n ,\end{equation}
\begin{equation}\label{2.4}f_n\rightarrow f \mbox{\ \ strongly \ \ in } L^m(Q).\end{equation}
According to the classical result(see \cite{29}), there exists at least a weak solution $u_n \in C([0, T]; L^{2}(\Omega))\cap L^{\vec{p}}(0,T;W_0^{1,\vec{p}}(\Omega))$ to
problem $(P_n)$.

In order to get some estimates of $u_n$, we need the following three results. The first one is the anisotropic version of Gagliardo-Nirenberg inequality(see \cite{32}) which is stated as follows.\vspace{2mm}

\noindent{\bf Proposition 3.1.} \emph{Let $\vec{\sigma}=(\sigma_1,\sigma_2,...,\sigma_N)$, where $\sigma_i>1$ for $i=1,2,...,N$ with $\bar{\sigma}<N$. If $u\in C\left(0,T;L^2(\Omega)\right)\cap L^{\vec{\sigma}}\left(0,T;W_{0}^{1,\vec{\sigma}}(\Omega)\right)$.
Then there exists a constant $C>0$ that depends only on $N$ and $\sigma_i$ such that
$$\iint_Q|u|^{\frac{2\bar{\sigma}}{N}+\bar{\sigma}}dxdt\leq C\left(\mathop {\mbox{ess\ sup}}_{t\in(0,T)}\int_\Omega|u|^2 dx\right)^{\frac{\bar{\sigma}}{N}}\prod_{i=1}^{N}\left(\iint_Q|D_iu|^{\sigma_i}dxdt\right)^{\frac{\bar{\sigma}}{\sigma_iN}}.$$}
%

%{\bf Proof.}
%Using H\"{o}lder's inequality and Lemma 2.1, it follows that
%\begin{equation*}\begin{split}
%\iint_Q|u|^{\frac{2\bar{\sigma}}{N}+\bar{\sigma}}dxdt
%&\leq\int_0^{T}\left(\int_\Omega|u|^{\frac{\bar{\sigma}N}{N-\bar{\sigma}}}dx\right)^{\frac{N-\bar{\sigma}}{N}}\left(\int_\Omega|u|^2dx\right)^{\frac{\bar{\sigma}}{N}}dt\\
%
%&\leq %C\int_0^T\prod_{i=1}^N\|D_iu\|_{L^{\sigma_i}(\Omega)}^{\frac{\bar{\sigma}}{N}}\left(\int_\Omega|u|^2dx\right)^{\frac{\bar{\sigma}}{N}}dt\\
%& \leq C\left(\mathop {\mbox{ess\ sup}}_{t\in(0,T)}\int_\Omega|u|^2 dx\right)
%^{\frac{\bar{\sigma}}{N}}\int_0^T\prod_{i=1}^{N}\|D_iu\|_{L^{\sigma_i}(\Omega)}^{\frac{\bar{\sigma}}{N}}dt\\
%& \leq C\left(\mathop {\mbox{ess\ sup}}_{t\in(0,T)}\int_\Omega|u|^2 dx\right)^{\frac{\bar{\sigma}}{N}}\prod_{i=1}^{N}\left(\iint_Q|D_iu|^{\sigma_i}dxdt\right)^{\frac{\bar{\sigma}}{\sigma_iN}},
%\end{split}\end{equation*}
%where $\bar{\sigma}^*=\frac{N\bar{\sigma}}{N-\bar{\sigma}}.$ \ \ $\Box$

The second one is the following embedding result of Marcinkiewicz space.\vspace{2mm}

\noindent{\bf Proposition 3.2.} \emph{Let $v$ be a measurable function belonging to ${\cal
M}^\mu(Q)$ for $\mu>0$, and suppose that there exist two
constants $\nu_i>\gamma\geq 0$ such that
\begin{align}\label{2.13}
\iint_{Q}|D_iT_k(v)|^{p_i} dxdt\leq
C(1+k)^{\gamma(p_i-1)} k^{\nu_i-\gamma(p_i-1)}, \  \ \forall k>0,
\end{align}
where $C$ is a
positive constant independent of $k$.
Then
$$D_iv\in {\cal M}^{\delta_i}(Q),$$ where
$\delta_i=\frac{p_i\mu}{\mu+\nu_i}$.}\vspace{2mm}

\noindent{\bf Proof.}  The proof is similar to the proof of Proposition 3.1 in \cite{d}(see also \cite{21}). For the sake of clarity and readability, we give the details
below.

For every
$k>0,$ we have
\begin{align}\label{2.14}
|\{|D_iv|>l\}|\leq |\{|v|>k\}|+|\{|D_iv|>l, |v|\leq
k\}|,\end{align}
where $l$ is a fixed positive number.

In view of (\ref{2.13}), it is easy to get
$$|\{|D_iv|>l, |v|\leq k\}|\leq\frac{1}{l^{p_i}}\iint_{Q}|D_iT_k(v)|^{p_i} dxdt\leq
C\frac{(1+k)^{\gamma(p_i-1)} k^{\nu_i-\gamma(p_i-1)}}{l^{p_i}}.$$
This implies
$$|\{|D_iv|>l, |v|\leq k\}|\leq
2^{\gamma(p_i-1)}C\frac{ k^{\nu_i}}{l^{p_i}},\ \ \mbox{for}\ k\geq1. $$
Noticing that $v\in {\cal
M}^\mu(Q)$, we have
$$|\{|v|>k\}|\leq \frac{\tilde{C}}{k^{\mu}},$$
where $\tilde{C}$ is a positive constant independent of $k$.

Combining the above two inequalities with (\ref{2.14}), we obtain
$$|\{|D_iv|>l\}|\leq \max\{2^{\gamma(p_i-1)}C, \tilde{C}\}\left(\frac{k^{\nu_i}}{l^{p_i}}+\frac{1}{k^{\mu}}\right),\ \ \mbox{for}\ k\geq1.$$
Minimizing with respect to $k$, namely taking $k={(\frac{\mu}{\nu_i})}^{\frac{1}{\mu+\nu_i}}l^{\frac{p_i}{\mu+\nu_i}}$, we get
$$|\{|D_iv|>l\}|\leq
\max\{2^{\gamma(p_i-1)}C, \tilde{C}\}\left[{(\frac{\mu}{\nu_i})}^{\frac{\nu_i}{\mu+\nu_i}}+1\right]\frac{ 1}{l^{\delta_i}},\ \ \mbox{for}\ l\geq{(\frac{\nu_i}{\mu})}^{\frac{1}{p_i}}.$$
Since $Q$ is bounded, the
above inequality obviously holds true for
$l<{(\frac{\nu_i}{\mu})}^{\frac{1}{p_i}}.$ Thus, the proof is finished. \ \ $\Box$\vspace{2mm}

Furthermore, by the proof of lemma 3.11 of \cite{35}, we easily get the following result.\vspace{2mm}

\noindent{\bf Proposition 3.3.} \emph{Let $r \geq1$ and $g :Q\rightarrow \mathbb{R}$ be a measurable function. If $$\sum_{\tilde{k}=0}^{\infty}\tilde{k}^{r-1}|\{|g|>\tilde{k}\}|<+\infty,$$
then $g\in L^r(Q)$ and $$\int_0^T\int_\xo|g|^rdxdt\leq 2^r|Q|+2^{r-1}\left[|Q|+\sum_{\tilde{k}=0}^{\infty}\tilde{k}^{r-1}|\{|g|>\tilde{k}\}|\right].$$}\vspace{2mm}

\noindent{\bf Lemma 3.1.} \emph{Suppose that (\ref{1.2})-(\ref{1.4}) hold true and $f\in
L^m(Q)$ with $m>\frac{N}{\bar{p}}+1$. If $u_n$ is a weak solution of problem $(P_n)$, then
 \begin{equation}\label{2.5}||u_n||_{L^\infty(Q)}\le C_1,
  \end{equation}
\begin{equation}\label{2.6} ||u_n||_{L^{\vec{p}}(0,T;W_0^{1,\vec{p}}(\Omega))}\le C_1. \end{equation}}\vspace{2mm}

\noindent{\bf Proof.} We only prove the case $\theta>0$, since the case $\theta=0$ is a well-known result.

For $k>0$, denote $G_k(s)=s-T_k(s), \forall
s\in \mathbb{R}, k>0$. Taking $G_k(u_n(x,
t))\chi_{(0,\tau)}(t)$ as a test function in problem $(P_n)$, and using (\ref{1.1}) and (\ref{2.3}) together with H\"{o}lder's inequality, we get for every $\tau\in (0,T]$
\begin{equation}\begin{split}\label{2.9} &\int_{A_k(\tau)}|G_k(u_n(\tau))|^2dx+
\int_0^\tau\int_{A_k(t)}\sum_{i=1}^{N}\frac{\alpha|D_iu_n|^{p_i}}{(1+|u_n|)^{\theta(p_i-1)}}dxdt\\
&\leq C\left(\int_0^\tau\int_{A_k(t)}|G_k(u_n)|^{m'}dxdt\right)^{\frac{1}{m'}},
\end{split}\end{equation}
where
\begin{equation}\label{0309}A_k(t)=\{x\in \Omega: |u_n(x,t)|>k\}.\end{equation}
Setting $$\sigma_i=\frac{p_i\bar{p}(N+2-N\theta)+N\theta p_i}{N\theta(p_i-\bar{p})+(N+2)\bar{p}},$$
then it is easy to check that $\sigma_i\leq p_i.$ Furthermore, we have $$\frac{\theta(p_i-1)\sigma_i}{p_i-\sigma_i}=\frac{(N+2)\bar{\sigma}}{N},\ \ \ \ \ \frac{p_i-\sigma_i}{p_i}=\frac{\theta N\sigma_i(p_i-1)}{(N+2)\bar{\sigma} p_i},$$ and
\begin{equation}\label{p3.7}
\bar{\sigma}=\frac{(N+2)\bar{p}-N\theta \bar{p}+N\theta}{N+2}.
\end{equation}

According to (\ref{2.9}) and H\"{o}lder's inequality, we get
\begin{equation*}
\begin{split}\label{2.10}
&\int_0^\tau\int_{A_k(t)}|D_iG_k(u_n)|^{\sigma_i}dxdt\\
&\leq\left(\int_0^\tau\int_{A_k(t)}\frac{|D_iu_n|^{p_i}}{(1+|u_n|)^{\theta(p_i-1)}}dxdt\right)^\frac{\sigma_i}{p_i}\left(\int_0^\tau\int_{A_k(t)}(1+|u_n|)^\frac{\theta (p_i-1)\sigma_i}{p_i-\sigma_i}dxdt\right)^\frac{p_i-\sigma_i}{p_i}\\
&\leq C\left(\int_0^\tau\int_{A_k(t)}|G_k(u_n)|^{m'}dxdt\right)^\frac{\sigma_i}{p_i m'}\left(\int_0^\tau\int_{A_k(t)}(k+|G_k(u_n)|)^\frac{(N+2)\bar{\sigma}}{N}dxdt\right)^\frac{N\theta\sigma_i(p_i-1)}{(N+2)\bar{\sigma} p_i}\\
&\leq C\left(\int_0^\tau\int_{A_k(t)}|G_k(u_n)|^{m'}dxdt\right)^\frac{\sigma_i}{p_i m'}\\
& \ \ \ \ \times\left[k^\frac{\theta\sigma_i(p_i-1)}{p_i}\left(\int_0^T|A_k(t)|dt\right)^{\frac{N\theta\sigma_i(p_i-1)}{(N+2)\bar{\sigma}p_i}}+\left(\int_0^\tau\int_{A_k(t)}|G_k(u_n)|^\frac{(N+2)\bar{\sigma}}{N}dxdt\right)^{\frac{N\theta\sigma_i(p_i-1)}{(N+2)\bar{\sigma}p_i}} \right].
\end{split}
\end{equation*}
From proposition 3.1, (\ref{2.9}) and
the above inequality,  it follows that
\begin{equation*}
\begin{split}
&\int_0^T\int_{A_k(t)} |G_k(u_n)|^{\frac{(N+2)\bar{\sigma}}{N}}dxdt\\
&\leq C\left(\mathop {\mbox{ess\ sup}}_{t\in(0,T)}\int_{A_k(t)}|G_k(u_n)|^2dx\right)^{\frac{\bar{\sigma}}{N}}\prod_{i=1}^{N}\left(\int_0^T\int_{A_k(t)}|D_iG_k(u_n)|^{\sigma_i} dxdt\right)^{\frac{\bar{\sigma}}{\sigma_iN}} \\
&\leq C\left(\int_0^T\int_{A_k(t)}
|G_k(u_n)|^{m'}dxdt\right)^{\frac{\bar{\sigma}(N+\bar{p})}{m'N\bar{p}}}\\
&\ \ \ \ \times
\left[k^{\frac{\theta\bar{\sigma}(\bar{p}-1)}{\bar{p}}}\left(\int_0^T|A_k(t)|dt\right)^{\frac{N\theta(\bar{p}-1)}{\bar{p}(N+2)}}
+\left(\int_0^T\int_{A_k(t)}|G_k(u_n)|^{\frac{(N+2)\bar{\sigma}}{N}}
 dxdt\right)^{\frac{N\theta(\bar{p}-1)}{\bar{p}(N+2)}}\right].
\end{split}
\end{equation*}
Recalling that $m>\frac{N}{\bar{p}}+1$ and $0\leq \theta <1+\frac{\bar{p}}{N(\bar{p}-1)},$ we get
$\frac{\bar{\sigma}(N+2)}{m'N}>1.$ Applying H\"{o}lder's inequality, we
have
\begin{equation*}\begin{split} &\int_0^T\int_{A_k(t)} |G_k(u_n)|^{\frac{(N+2)\bar{\sigma}}{N}}dxdt\\
 &\le C\left(\int_0^T\int_{A_k(t)}
|G_k(u_n)|^{\frac{(N+2)\bar{\sigma}}{N}}dxdt\right)^{\frac{N+\bar{p}}{\bar{p}(N+2)}}\left(\int_0^T|A_k(t)|dt\right)^{\frac{\bar{\sigma}(N+\bar{p})(m-1)}{\bar{p}Nm}-\frac{N+\bar{p}}{\bar{p}(N+2)}}\\
& \ \ \ \ \times
\left[k^{\frac{\theta\bar{\sigma}(\bar{p}-1)}{\bar{p}}}\left(\int_0^T|A_k(t)|dt\right)^{\frac{N\theta(\bar{p}-1)}{\bar{p}(N+2)}}
+\left(\int_0^T\int_{A_k(t)}|G_k(u_n)|^{\frac{(N+2)\bar{\sigma}}{N}}
 dxdt\right)^{\frac{N\theta(\bar{p}-1)}{\bar{p}(N+2)}}\right],\\
\end{split}\end{equation*}
which implies that
\begin{equation*}\begin{split}\label{2.15}
 &\int_0^T\int_{A_k(t)}
 |G_k(u_n)|^{\frac{(N+2)\bar{\sigma}}{N}}dxdt\\
 &\leq
 C\left(\int_0^T|A_k(t)|dt\right)^{\frac{\bar{\sigma}(N+\bar{p})(m-1)(N+2)}{Nm[N(\bar{p}-1)+\bar{p}]}-\frac{N+\bar{p}}{N(\bar{p}-1)+\bar{p}}}\\
& \ \ \ \ \times\left[k^{\frac{\theta(\bar{p}-1)\bar{\sigma}(N+2)}{N(\bar{p}-1)+\bar{p}}}\left(\int_0^T|A_k(t)|dt\right)^{\frac{N\theta(\bar{p}-1)}{N(\bar{p}-1)+\bar{p}}}
 +\left(\int_0^T\int_{A_k(t)}|G_k(u_n)|^{\frac{(N+2)\bar{\sigma}}{N}}
  dxdt\right)^{\frac{N\theta(\bar{p}-1)}{N(\bar{p}-1)+\bar{p}}}\right].
\end{split}\end{equation*}
Therefore, by Young's inequality and the fact that $\frac{N\theta(\bar{p}-1)}{N(\bar{p}-1)+\bar{p}}<1,$ we obtain
\begin{align}\label{2.17}
&\int_0^T\int_{A_k(t)}|G_k(u_n)|^{\frac{(N+2)\bar{\sigma}}{N}}dxdt\notag\\
&\leq Ck^{\frac{\theta(\bar{p}-1)\bar{\sigma}(N+2)}{N(\bar{p}-1)+\bar{p}}}\left(\int_0^T|A_k(t)|dt\right)^{\lambda}\notag\\
&\ \ \ \ \ +C\left(\int_0^T|A_k(t)|dt\right)^{\left[\frac{\bar{\sigma}(N+\bar{p})(m-1)(N+2)}{Nm[N(\bar{p}-1)+\bar{p}]}
-\frac{N+\bar{p}}{N(\bar{p}-1)+\bar{p}}\right]\frac{N(\bar{p}-1)+\bar{p}}{N(\bar{p}-1)-N\theta(\bar{p}-1)+\bar{p}}},
\end{align}
where
 $$\lambda=\frac{\bar{\sigma}(N+\bar{p})(m-1)(N+2)}{Nm[N(\bar{p}-1)+\bar{p}]}
+\frac{N\theta(\bar{p}-1)-N-\bar{p}}{N(\bar{p}-1)+\bar{p}}.$$
Since $m>\frac{N}{\bar{p}}+1$, it is easy to check that \begin{align}\left[\frac{\bar{\sigma}(N+\bar{p})(m-1)(N+2)}{Nm[N(\bar{p}-1)+\bar{p}]}-\frac{N+\bar{p}}{N(\bar{p}-1)+
\bar{p}}\right]\frac{N(\bar{p}-1)+\bar{p}}{N(\bar{p}-1)-N\theta(\bar{p}-1)+\bar{p}}>\lambda>1\notag,\end{align}
where we have used (\ref{p3.7}).

Hence, we get for $k\geq1$
\begin{equation*}\begin{split}\label{2.18}
 &\int_0^T\int_{A_k(t)}
 |G_k(u_n)|^{\frac{(N+2)\bar{\sigma}}{N}}dxdt\\
&\leq Ck^{\frac{\theta(\bar{p}-1)\bar{\sigma}(N+2)}{N(\bar{p}-1)+\bar{p}}}
 \left(\int_0^T|A_k(t)|dt\right)^{\lambda} \\
&\ \ \ \
\mbox{+}C|Q|^{\left[\frac{\bar{\sigma}(N+\bar{p})(m-1)(N+2)}{Nm[N(p-1)+p]}-\frac{N+\bar{p}}{N(\bar{p}-1)+\bar{p}}\right]
\frac{N(\bar{p}-1+\bar{p})}{N(\bar{p}-1)-N\theta(\bar{p}-1)+\bar{p}}-\lambda}\left(\int_0^T|A_k(t)|dt\right)^{\lambda}\\
&\leq
 Ck^{\frac{\theta(\bar{p}-1)\bar{\sigma}(N+2)}{N(\bar{p}-1)+\bar{p}}}\left(\int_0^T|A_k(t)|dt\right)^{\lambda},\\
 \end{split}\end{equation*}
which yields
\begin{equation}\label{2.21}
\varphi(h)\leq\frac{Ck^{\frac{\theta(\bar{p}-1)\bar{\sigma}(N+2)}{N(\bar{p}-1)+\bar{p}}}
\varphi(k)^{\lambda}}{(h-k)^{\frac{(N+2)\bar{\sigma}}{N}}},
\ \ \ \forall h>k\geq 1,\end{equation}
where
\begin{equation}\notag\varphi(k)=\int_0^T|A_k(t)|dt. \end{equation}
By Lemma A.2 in
\cite{21}, one may conclude that there exists a  constant $k^*>0$ independent of $n$ satisfying
\begin{equation}\label{2.22}\varphi(k^*)=0, \end{equation}
from which, we get (\ref{2.5}). Taking $u_n$ as a test function in problem $(P_n)$ and using (\ref{2.5}),
we easily get (\ref{2.6}).\ \ $\Box$\vspace{2mm}

\noindent{\bf Lemma 3.2.} \emph{Assume that (\ref{1.2})-(\ref{1.4}) hold true and $f\in
L^m(Q)$ with $m=\frac{N}{\bar{p}}+1$. If $u_n$ is a weak solution of problem $(P_n)$, then
\begin{equation}\label{2.23}||u_n||_{L^{\vec{p}}(0,T;W_0^{1,\vec{p}}(\Omega))}\leq C_2, \end{equation}
\begin{equation}\label{2.24}\|e^{\lambda |u_n|^{1-\frac{N\theta(\bar{p}-1)}{N(\bar{p}-1)+\bar{p}}}}\|_{L^1(Q)}\leq C_2.\end{equation}}\vspace{2mm}

\noindent{\bf Proof. }In the same way as in the proof of (\ref{2.17}), we get
\begin{align}
&\int_0^T\int_{A_k(t)}|G_k(u_n)|^{\frac{(N+2)\bar{\sigma}}{N}}dxdt
\leq Ck^{\frac{\theta(\bar{p}-1)\bar{\sigma}(N+2)}{N(\bar{p}-1)+\bar{p}}}\int_0^T|A_k(t)|dt+C\int_0^T|A_k(t)|dt,\notag
\end{align}
and as a result, one obtains
\begin{equation*}\varphi(h)\leq\frac{Ck^{\frac{\theta(\bar{p}-1)
\bar{\sigma}(N+2)}{N(\bar{p}-1)+\bar{p}}}\varphi(k)}{(h-k)^{\frac{(N+2)\bar{\sigma}}{N}}},
\ \ \ \forall h>k\geq 1.\end{equation*}
This, combined with lemma 3 of \cite{kov05}(here $\tau_1=\frac{\theta(\bar{p}-1)
\bar{\sigma}(N+2)}{N(\bar{p}-1)+\bar{p}}$ and $\tau_2=\frac{(N+2)\bar{\sigma}}{N}$), gives
\begin{align}
&\int_{0}^T|\{|u_n|>k\}|dt\notag\\
&\leq e^{s_0+1}e^{-\left({\frac{k-1}{\lambda_0}}\right)^{1-\frac{N\theta(\bar{p}-1)}{N(\bar{p}-1)+\bar{p}}}}\int_{0}^T|\{|u_n|>1\}|dt ,
\ \ \ \forall k\geq 1+\lambda_0s_0^{\frac{\tau_2}{\tau_2-\tau_1}},\notag\end{align}
where $\lambda_0=\mbox{max}\{1,\left[(Ce)^{\frac{1}{\tau_2}}2^{\frac{\tau_1}{\tau_2}}
\frac{\tau_2-\tau_1}{\tau_2}\right]^{\frac{\tau_2}{\tau_2-\tau_1}}\}$, and $s_0\in \mathbb{N}$ with $s_0>\frac{2^{\frac{\tau_2-\tau_1}{\tau_2}}}{2^{\frac{\tau_2-\tau_1}{\tau_2}}-1}$.

Note that $k\geq 1+\lambda_0s_0^{\frac{\tau_2}{\tau_2-\tau_1}}\geq2$ implies that $k-1\geq \frac{k}{2}.$ As a result, the above inequality shows that
\begin{equation*}\int_{0}^T|\{|u_n|>k\}|dt\leq Ce^{-2\lambda{k}^{1-\frac{N\theta(\bar{p}-1)}{N(\bar{p}-1)+\bar{p}}}},
\end{equation*}
where $2\lambda=(\frac{1}{2\lambda_0})^{1-\frac{N\theta(\bar{p}-1)}{N(\bar{p}-1)+\bar{p}}}.$
Thus
\begin{align}\label{316}&\int_0^T\left|\left\{e^{\lambda |u_n|^{1-\frac{N\theta(\bar{p}-1)}{N(\bar{p}-1)+\bar{p}}}}>e^{\lambda k^{1-\frac{N\theta(\bar{p}-1)}{N(\bar{p}-1)+\bar{p}}}}\right\}\right|dt\notag\\
&=\int_{0}^T|\{|u_n|>k\}|dt\leq C\left(e^{\lambda{k}^{1-\frac{N\theta(\bar{p}-1)}{N(\bar{p}-1)+\bar{p}}}}\right)^{-2},\ \ \forall k\geq 1+\lambda_0s_0^{\frac{\tau_2}{\tau_2-\tau_1}}.
\end{align}
Now let us take $k=\left(\frac{\mbox{ln} \tilde{k}}{\lambda}\right)^{\frac{1}{1-\frac{N\theta(\bar{p}-1)}{N(\bar{p}-1)+\bar{p}}}}
$ with $\tilde{k}\in \mathbb{N}$ satisfying
$$\tilde{k}\geq \mbox{max}\{e^{\lambda{2}^{1-\frac{N\theta(\bar{p}-1)}{N(\bar{p}-1)+\bar{p}}}},
e^{\lambda{(1+\lambda_0s_0^{\frac{\tau_2}{\tau_2-\tau_1}})}^{1-\frac{N\theta(\bar{p}-1)}{N(\bar{p}-1)+\bar{p}}}}\}=\bar{k}.$$
Then, by (\ref{316}), it follows that
\begin{equation*}
\int_0^T\left|\left\{e^{\lambda |u_n|^{1-\frac{N\theta(\bar{p}-1)}{N(\bar{p}-1)+\bar{p}}}}>\tilde{k}\right\}\right|dt\leq \frac{|Q|e}{\tilde{k}^2},\ \ \forall \tilde{k}\geq \bar{k},
\end{equation*}
which yields
$$\sum_{\tilde{k}\in\mathbb{N},\tilde{k}\geq \bar{k}}\int_0^T\left|\left\{e^{\lambda |u_n|^{1-\frac{N\theta(\bar{p}-1)}{N(\bar{p}-1)+\bar{p}}}}>\tilde{k}\right\}\right|dt\leq |Q|e\sum_{\tilde{k}=1}^\infty\frac{1}{\tilde{k}^2}<+\infty.$$
This, combined with proposition 3.3, gives (\ref{2.24}).

 Let us take $\psi(u_n)=[(1+|u_n|)^l-1]\mbox {sign}(u_n)\chi_{(0,\tau)}(t)$ as a test function of problem $(P_n)$, where $l>1$ will be chosen later. It follows from (\ref{1.2}) and integrating by parts that
\begin{align}\label{p0314}
 &\int_\Omega \Psi(u_n(x,\tau))dx+\alpha l\sum_{i=1}^{N}\iint_{Q}
\frac{|D_iu_n|^{p_i}}{(1+|u_n|)^{\theta(p_i-1)+1-l}}dxdt\notag\\
&\leq\int_0^\tau\int_\Omega |f_n|(1+|u_n|^l)dxdt,
\end{align}
where $\Psi(s)=\int_0^{s}\psi(\xi)d\xi.$
Observing that
\begin{equation*}\Psi(s)\geq \frac{1}{2l+1}|s|^{l+1}-C, \ \ \ \forall s\in \mathbb{R},
\end{equation*}
and taking the supremum of both terms for $\tau\in(0, T)$ in (\ref{p0314}), we obtain
\begin{equation}\label{2.27}
\begin{split}
&\frac{1}{2l+1}\mathop {\mbox{ess\ sup}}_{0\leq t\leq T}\int_\Omega|u_n(x, t)|^{l+1}dx+\alpha l\sum_{i=1}^{N}\iint_{Q}
\frac{|D_iu_n|^{p_i}}{(1+|u_n|)^{\theta(p_i-1)+1-l}}dxdt\\
&\leq\|f\|_{L^m(Q)}\left(\iint_Q(1+|u_n|)^{lm^{'}}dxdt\right)^{\frac{1}{m^{'}}}+C.
\end{split}
\end{equation}
Taking $l>1+\theta(\mbox{max}\{p_1,p_2,...,p_N\}-1)$, using (\ref{2.24}) and (\ref{2.27}), we have
\begin{equation*}
\begin{split}
\sum_{i=1}^{N}\iint_{Q}|D_iu_n|^{p_i}dxdt&\leq\sum_{i=1}^{N}\iint_{Q}|D_iu_n|^{p_i}(1+|u_n|)^{l-1-\theta(p_i-1)}dxdt\\
&\leq C\left(\iint_{Q}(1+|u_n|)^{lm'}dxdt\right)^{\frac{1}{m'}}+C\\
&\leq C,
\end{split}
\end{equation*}
from which, we get the desired result (\ref{2.23}). \ \ $\Box$\vspace{2mm}

\noindent{\bf Lemma 3.3.} \emph{Assume that (\ref{1.2})-(\ref{1.4}) hold true and $f\in
L^m(Q)$ with $m$ satisfying condition (\ref{1.03}). If $u_n$ is a weak solution of problem $(P_n)$, then
\begin{equation}\label{2.41}||u_n||_{L^{\vec{p}}(0,T;W_0^{1,\vec{p}}(\Omega))}\leq C_{3}, \end{equation}
\begin{equation}\label{2.42}||u_n||_{L^{\tilde{r}}(Q)}\leq C_{3},\end{equation}
where $\tilde{r}$ is defined as in (\ref{1.04}).}\vspace{2mm}

\noindent{\bf Proof.} We also use $\psi=[(1+|u_n|)^l-1]\mbox {sign}(u_n)\chi_{(0,\tau)}(t)$ as a test function of problem $(P_n)$, and then obtain (\ref{2.27}), where $l>1$ will be chosen later.
Set
$$q_i=\frac{[\bar{p}(N-N\theta+l+1)+N(\theta+l-1)]p_i}{N\theta p_i+N(\lambda-1)p_i+\bar{p}(N-N\theta+l+1)},$$
where $1<\lambda<l.$ A careful calculation shows that
\begin{equation}\label{p3.14}
1<q_i<p_i\ \ \ \mbox{and}\ \ \  \frac{q_i[\theta(p_i-1)+1-l]}{p_i-q_i}+(\lambda-1)\frac{p_iq_i}{p_i-q_i}
=\frac{N+\frac{l+1}{\lambda}}{N}\lambda\bar{q},
\end{equation}
where
\begin{equation}\label{p0317}
\bar{q}=\frac{\bar{p}(N-N\theta+l+1)+N(\theta+l-1)}{N\lambda+l+1}.
\end{equation}
Combining (\ref{2.27}), (\ref{p3.14}) with H\"{o}lder's inequality, we deduce
\begin{equation*}
\begin{split}
&\iint_Q |D_i|u_n|^\lambda|^{q_i}dxdt=\lambda^{q_i}\iint_Q|u_n|^{(\lambda-1)q_i}|D_i u_n|^{q_i}dxdt\\
&\leq \lambda^{q_i}\left(\iint_Q\frac{|D_i u_n|^{p_i}}{(1+|u_n|)^{\theta(p_i-1)+1-l}}dxdt\right)^{\frac{q_i}{p_i}}\left(\iint_Q (1+|u_n|)^{[\frac{\theta(p_i-1)+1-l}{p_i}+(\lambda-1)]\frac{q_i p_i}{p_i-q_i}} dxdt\right)^{\frac{p_i-q_i}{p_i}}\\
&\leq C\left(\left(\iint_Q(1+|u_n|)^{lm'}dxdt\right)^{\frac{q_i}{m'p_i}}+1\right)\left(\iint_Q(1+|u_n|)^{\frac{N+\frac{l+1}{\lambda}}{N}\lambda \bar{q}}dxdt\right)^{\frac{p_i-q_i}{p_i}}.
\end{split}
\end{equation*}
Now let
$$lm'=\frac{N\lambda+l+1}{N}\bar{q},$$
that is
\begin{equation}\label{p0318}
l=\frac{(m-1)[N(1-\theta)(\bar{p}-1)+\bar{p}]}{N+\bar{p}-\bar{p}m}.
\end{equation}
Hence, applying Young's inequality in the above estimate gives
\begin{equation}
\begin{split}\label{2.33}
&\iint_Q |D_i|u_n|^\lambda|^{q_i}dxdt\\
&\leq C\left[\left(\iint_Q(1+|u_n|)^{\frac{N\lambda+l+1}{N}\bar{q}}dxdt\right)
^{\frac{q_i}{m'p_i}+\frac{p_i-q_i}{p_i}}+1\right].
\end{split}
\end{equation}
Denoting
\begin{equation}\label{p3.17}
r=\frac{N\lambda+l+1}{N\lambda}\bar{q},\end{equation}
we get the following interpolation inequality for $L^r$-norms of $|u_n|^{\lambda}$
\begin{equation}\label{2.34}
\||u_n|^{\lambda}\|_{L^r(\Omega)}\leq\||u_n|^{\lambda}\|_{L^{\bar{q}^*}(\Omega)}
^{\delta}\||u_n|^{\lambda}\|_{L^{\frac{l+1}{\lambda}}(\Omega)}^{1-\delta},
\end{equation}
with $\delta$ satisfying
\begin{equation*}\frac{1}{r}=\frac{\delta}{\bar{q}^*}+\frac{\lambda(1-\delta)}{l+1}.\end{equation*}
This, combined with (\ref{p3.17}),  gives
\begin{equation}\label{p3.19}
\delta=\frac{N\lambda}{N\lambda+l+1},\ \ \ \ \delta r=\bar{q}.
\end{equation}
%\begin{equation}\label{2.35}
%r=\frac{N\lambda+l+1}{N\lambda}\bar{q}=\frac{\bar{p}(l+1)+N[l+(1-\theta)(\bar{p}-1)]}{N\lambda}.
%\end{equation}
%Hence and $\delta r=\bar{q}.$
Exploiting (\ref{2.27}), (\ref{2.34}) and  (\ref{p3.19}), and using Lemma 2.1, we get
\begin{equation*}
\begin{split}
\int_{0}^{T}\left\||u_n|^{\lambda}\right\|_{L^r(\Omega)}^{r}dt
&\leq\int_{0}^{T}\left\||u_n|^{\lambda}\right\|_{{L^{\bar{q}^*}}(\Omega)}^{\delta r}
\left\||u_n|^{\lambda}\right\|_{{L^{\frac{l+1}{\lambda}}}(\Omega)}^{(1-\delta)r}dt\\
&\leq C\left[\left(\iint_Q(1+|u_n|)^{lm^{'}}dxdt\right)^{\frac{\lambda(1-\delta)r}{m^{'}(l+1)}}
+1\right]\int_{0}^{T}\prod_{i=1}^{N}\left\|D_i|u_n|^{\lambda}\right\|_{L^{q_i}(\Omega)}^{\frac{\bar{q}}{N}}dt.
\end{split}
\end{equation*}
Observing that $\sum\limits_{i=1}^{N}\frac{\bar{q}}{N q_i}=1$, using general H\"{o}lder's inequality and Young's inequality, then it follows from the above estimate that
\begin{align}
\label{2.36}
&\int_{0}^{T}\left\||u_n|^{\lambda}\right\|_{L^r(\Omega)}^{r}dt\notag\\
&\leq C\left[\left(\iint_Q(1+|u_n|)^{lm^{'}}dxdt\right)^{\frac{\lambda(1-\delta)r}{m^{'}(l+1)}}
+1\right]\prod_{i=1}^{N}\left(\iint_Q\left|D_i|u_n|^{\lambda}\right|^{q_i}dxdt\right)^{\frac{\bar{q}}{Nq_i}}\notag\\
&\leq C\left[\left(\iint_Q(1+|u_n|)^{lm^{'}}dxdt\right)^{\frac{\lambda(1-\delta)r}{m^{'}(l+1)}}
+1\right]\prod_{i=1}^{N}\left[\left(\iint_Q(1+|u_n|)^{lm'} dxdt\right)^{\frac{q_i}{m'p_i}+\frac{p_i-q_i}{p_i}}+1\right]^{\frac{\bar{q}}{Nq_i}}\notag\\
&\leq C\left(\iint_Q(1+|u_n|)^{lm'} dxdt\right)^{\frac{\lambda(1-\delta)r}{m'(l+1)}+\sum\limits_{i=1}^{N}\left[{\frac{q_i}{m'p_i}+\frac{p_i-q_i}{p_i}}\right]\frac{\bar{q}}{N q_i}}+C.
\end{align}
Let $\tilde{r}=r\lambda$, then it is apparent that $\tilde{r}=lm'=\frac{mN[(1-\theta)(\bar{p}-1)]+m\bar{p}}{N-m\bar{p}+\bar{p}}$. Therefore, by (\ref{2.36}), it follows that
\begin{equation*}
 \iint_Q|u_n|^{\tilde{r}}dxdt\leq C\left(\iint_Q|u_n|^{\tilde{r}} dxdt\right)^{\tilde{h}}+C,
\end{equation*}
where
\begin{equation*}
\begin{split}
\tilde{h}=\frac{l(1-\delta)}{l+1}+\frac{(m-1)\bar{q}}{m\bar{p}}+1-\frac{\bar{q}}{\bar{p}}.
\end{split}
\end{equation*}
By (\ref{p0317}), (\ref{p0318}), (\ref{p3.19}) and the condition $m<\frac{N}{\bar{p}}+1$ , one can check that $\tilde{h}<1.$ Therefore, by Young's inequality we get (\ref{2.42}) immediately.

Note that condition (\ref{1.03}) gives
$l\geq 1+\theta\left(\max\limits_{1\leq i\leq N}p_i-1\right)$. This, combined with (\ref{2.27}) and (\ref{2.42}), gives (\ref{2.41}). The proof is completed. \ \ $\Box$ \vspace{2mm}

\noindent{\bf Lemma 3.4.} \emph{Assume that (\ref{1.2})-(\ref{1.4}) hold true and $f\in
L^m(Q)$ with $m$ satisfying condition (\ref{1.5}). If $u_n$ is a weak solution of problem $(P_n)$, then
\begin{equation}\label{2.48}||u_n||_{L^{\widetilde{r}}(Q)}\le C_{4},\end{equation}
\begin{equation} \label{2.47}||u_n||_{L^{\vec{q}}(0, T; W_0^{1, \vec{q}}(\Omega))}\le C_{4},\end{equation}
%\begin{equation}\label{2.49}||u_n^\prime||_{L^{\vec{q'}}(0,T;W_0^{-1,\vec{q'}}(\Omega))+L^m(Q)}\le C_{4},\end{equation}
with $\vec{q}=(q_1,q_2,...,q_N)$, where $ q_i$ and $\tilde{r}$ are defined as in (\ref{1.04}) and (\ref{2.08}) respectively.}\vspace{2mm}

\noindent{\bf Proof.} In the same way as in
the proof of Lemma 3.3, taking $\psi=[(1+|u_n|)^l-1]\mbox {sign}(u_n)\chi_{(0,\tau)}(t)$ as a test function of problem $(P_n)$, we get (\ref{2.27}), where $l$ will be fixed later.

Let us set
\begin{equation}\label{b}
q_i =\frac{p_i[\bar{p}(N+l+1)-N\theta\bar{p}+N(\theta+l-1)]}{N[\theta(p_i-1)-l+1]+\bar{p}(N+l+1)-N\theta\bar{p}+N(\theta+l-1)},
\end{equation}
then $q_i<p_i$ and
\begin{equation}\label{54}\frac{[\theta(p_i-1)-l+1]q_i}{p_i-q_i}=\frac{N+l+1}{N}\bar{q},\end{equation}
where
\begin{equation}\label{2.54}
\bar{q}=\frac{\bar{p}(N+l+1)-N\theta\bar{p}+N(\theta+l-1)}{N+l+1}.
\end{equation}
Thanks to (\ref{2.27}) and (\ref{54}), and using H\"{o}lder's inequality, we obtain
\begin{equation*}\begin{split}
  &\iint_Q|D_iu_n|^{q_i} dxdt\\
 &\le\left(\iint_Q \frac{|D_iu_n|^{q_i}}{(1+|u_n|)^{\theta(p_i-1)+1-l}}dxdt\right)^{\frac{q_i}{p_i}}
 \left(\iint_Q(1+|u_n|)^{\frac{[\theta(p_i-1)+1-l]q_i}{p_i-q_i}} dxdt\right)^{\frac{p_i-q_i}{p_i}}\\
&\leq\left[C\left(\iint_Q
|u_n|^{lm'}dxdt\right)^{\frac{1}{m'}}+C\right]^{\frac{q_i}{p_i}}\left(\iint_Q(1+|u_n|)^{{\frac{N+l+1}{N}}\bar{q}} dxdt\right)^{\frac{p_i-q_i}{p_i}}.
\end{split}\end{equation*}
This, combined with Young's inequality, shows
\begin{equation}\label{2.57} \iint_Q|D_iu_n|^{q_i} dxdt
 \le C\left(\iint_Q
|u_n|^{lm'}dxdt\right)^{\frac{q_i}{m'p_i}+\frac{p_i-q_i}{p_i}}+C,
\end{equation}
where we have set
\begin{equation}\label{2.56}lm'=\frac{(N+1+l)\bar{q}}{N}.
\end{equation}
From (\ref{b}), (\ref{2.54}) and (\ref{2.56}), it follows that
\begin{equation*}
\begin{split}
 &l=\frac{(m-1)[N(1-\theta)(\bar{p}-1)+\bar{p}]}{N+\bar{p}-\bar{p}m},\\
 &\bar{q}=\frac{m[N(1-\theta)(\bar{p}-1)+\bar{p}]}{N+1-[1+\theta(\bar{p}-1)](m-1)},
\end{split}
\end{equation*}
and
\begin{equation*}
\begin{split}\label{a}
q_i=\frac{m  p_i[N(1-\theta)(\bar{p}-1)+\bar{p}]}{[N-\bar{p}(m-1)][\theta(p_i-1)+1]+N(\bar{p}-1)(1-\theta)+\bar{p}}.
\end{split}
\end{equation*}
Now, setting
\begin{equation}\label{2.59}
\tilde{r}=\frac{N+l+1}{N}\bar{q}=\frac{\bar{p}(l+1)+N[l+(1-\theta)(\bar{p}-1)]}{N},
\end{equation}
we get the interpolation inequality for $L^{\tilde{r}}$-norms of $u_n$
\begin{equation}\label{2.58}
\|u_n\|_{L^{\widetilde{r}}(\Omega)}\leq\|u_n\|_{L^{\bar{q}^*}(\Omega)}^{\delta}\|u_n\|_{L^{l+1}(\Omega)}^{1-\delta},
\end{equation}
with $0<\delta<1$ such that
$$\frac{1}{\tilde{r}}=\frac{\delta}{\bar{q}^*}+\frac{1-\delta}{l+1}.$$
This, together with (\ref{2.59}), indicates
\begin{equation}
\label{p334}
\delta=\frac{N^2(l+1)-N\bar{q}(N+l+1)}{[(N-\bar{q})(l+1)-N\bar{q}](N+l+1)}=\frac{N}{N+l+1}.
\end{equation}
In view of (\ref{2.59})-(\ref{p334}) and Lemma 2.1, and reasoning exactly as in the proof of (\ref{2.36}), we get
\begin{equation*}
\begin{split}
\int_0^T\|u_n\|_{L^{\tilde{r}}(\Omega)}^{\tilde{r}}dt
&\leq\int_0^T\|u_n\|_{L^{\bar{q}^*}(\Omega)}^{\delta\tilde{r}}\|u_n\|_{L^{l+1}(\Omega)}^{(1-\delta)\tilde{r}}dt\\
&\leq C\left[\left(\iint_Q|u_n|^{lm'}dxdt\right)^{\frac{(1-\delta)\tilde{r}}{m'(l+1)}}
+1\right]\int_0^T\|u_n\|_{L^{\bar{q}^*}(\Omega)}^{\bar{q}}dt\\
&\leq C\left[\left(\iint_Q|u_n|^{lm'}dxdt\right)^{\frac{(1-\delta)\tilde{r}}{m'(l+1)}}
+1\right]\int_0^T\prod_{i=1}^{N}\|D_iu_n\|_{L^{q_i}(\Omega)}^{\frac{\bar{q}}{N}}dt\notag\\
&\leq C\left(\iint_Q|u_n|^{lm'}dxdt\right)^{\frac{(1-\delta)\tilde{r}}{m'(l+1)}+\sum\limits_{i=1}^N
\left[\frac{q_i}{m'p_i}+\frac{p_i-q_i}{p_i}\right]\frac{\bar{q}}{Nq_i}}+C.
\end{split}
\end{equation*}
Note that $lm'=\tilde{r}$ and $\frac{(1-\delta)\tilde{r}}{m'(l+1)}+\sum\limits_{i=1}^N
\left[\frac{q_i}{m'p_i}+\frac{p_i-q_i}{p_i}\right]\frac{\bar{q}}{Nq_i}<1$. Using Young's inequality, we can conclude that
(\ref{2.48}) holds true.

By (\ref{2.48}), (\ref{2.57}) and using the fact $lm'=\tilde{r}$, we get the result (\ref{2.47}).
The proof is completed.\ \ $\Box$ \vspace{2mm}

\noindent{\bf Lemma 3.5.} \emph{Assume that (\ref{1.2})-(\ref{1.4}) hold true and $f\in
L^m(Q)$ with $m$ satisfying condition (\ref{1.9}). If $u_n$ is a weak solution of problem $(P_n)$, then
\begin{equation} \label{3.59}||u_n||_{L^{\infty}(0, T; L^{1}(\Omega))}\leq C_{5}, \end{equation}
\begin{equation}\label{3.60}\|D_iT_k(u_n)\|_{L^{p_i}(\Omega)}\leq C_{5}(1+k)^{\frac{\theta(p_i-1)+1}{p_i}},\end{equation}
\begin{equation}\label{3.61}|\left\{|u_n|>k\right\}|\leq \frac{C_{5}}{k^{\tilde{r}}},  \end{equation}
\begin{equation}\label{3.62}|\left\{|D_iu_n|>k\right\}|\leq \frac{C_{5}}{k^{q_i}},\end{equation}
where $ q_i$ and $\tilde{r}$ are defined as in (\ref{1.04}) and (\ref{2.08}) respectively.}\vspace{2mm}

\noindent{\bf Proof.} The proof is divided into two steps.

{\bf Step 1.} In this step, we give the proof of the estimates of (\ref{3.61}) and (\ref{3.62}).

For all $\tau\in (0,T]$, choosing $T_k(u_n(x,t))\chi_{(0,\tau)}(t)$ as a test function for problem $(P_n)$, and using (1.2) and H\"{o}lder's inequality, we have
\begin{equation}
\begin{split}\label{3.64}
&\mathop {\mbox {ess\ sup}}_{0\leq t\leq T}\int_\Omega|T_k(u_n(x,\tau))|^2dx+\iint_Q\sum_{i=1}^{N}\frac{|D_iT_k(u_n)|^{p_i}}{(1+|u_n|)^{\theta(p_i-1)}}dxdt\\
&\leq C\left(\iint_Q|T_k(u_n)|^{m'}dxdt\right)^{\frac{1}{m'}}.
\end{split}
\end{equation}
Hence
\begin{equation}
\begin{split}\label{3.65}
&\iint_Q|D_iT_k(u_n)|^{p_i}dxdt\leq C (1+k)^{\theta(p_i-1)}\left(\iint_Q|T_k(u_n)|^{m'}dxdt\right)^{\frac{1}{m'}}.
\end{split}
\end{equation}
 The rest of the proof is divided into three cases.

{\bf case (i):} $m>\frac{\bar{p}(N+2)}{\bar{p}(N+2)-N}.$

In this case, we have $m'<\frac{\bar{p}(N+2)}{N}$. Hence we can choose $\rho_i$ such that
\begin{equation}\notag
p_i\geq\rho_i\geq\frac{m'p_i}{\theta(p_i-1)+m'},\ \ \ \frac{\bar{\rho}(N+2)}{N}=m',
\end{equation}
i.e.
\begin{equation}\label{3.67}
\frac{\theta(p_i-1)\rho_i}{p_i-\rho_i}\geq m',\ \ \ \ \bar{\rho}=\frac{Nm}{(N+2)(m-1)}.
\end{equation}
Using proposition 3.1, H\"{o}lder's inequality and (\ref{3.67}), we obtain that for $k\geq1$
\begin{equation*}
\begin{split}\label{3.69}
&\iint_Q|T_k(u_n)|^{\frac{(N+2)\bar{\rho}}{N}}dxdt\\
&\leq C\left(\mathop {\mbox {ess\ sup}}_{0\leq t\leq T}\int_\Omega|T_k(u_n)|^2dx\right)^{\frac{\bar{\rho}}{N}}\prod_{i=1}^{N}\left[\iint_Q|D_i T_k(u_n)|^{\rho_i}dxdt\right]^{\frac{\bar{\rho}}{\rho_iN}}\\
&\leq C\left(\iint_Q|T_k(u_n)|^{m'}dxdt\right)^{\frac{\bar{\rho}}{Nm'}}\prod_{i=1}^{N} \left(\iint_Q\frac{|D_i T_k(u_n)|^{p_i}}{(1+|T_k(u_n)|)^{\theta(p_i-1)}}dxdt\right)^{\frac{\bar{\rho}}{p_iN}}\\
& \ \ \ \ \times\left[\iint_Q(1+|T_k(u_n)|)^{\frac{\theta(p_i-1)\rho_i}{p_i-\rho_i}-m'}(1+|T_k(u_n)|)^{m'}dxdt\right]^{\frac{p_i-\rho_i}{p_i}\frac{\bar{\rho}}{\rho_iN}}\\
&\leq C\left(\iint_Q|T_k(u_n)|^{m'}dxdt\right)^{\frac{\bar{\rho}}{Nm'}}\\
& \ \ \ \ \times\prod_{i=1}^{N}(2k)^{\left[\frac{\theta(p_i-1)\rho_i}{p_i-\rho_i}-m'\right]\frac{p_i-\rho_i}{p_i}\frac{\bar{\rho}}{\rho_iN}}
\left[\iint_Q(1+|T_k(u_n)|)^{m'}dxdt\right]^{\frac{\bar{\rho}}{m'p_iN}+\frac{p_i-\rho_i}{p_i}\frac{\bar{\rho}}{\rho_iN}}.
\end{split}
\end{equation*}
If $\iint_Q|T_k(u_n)|^{m'}dxdt\geq 1$, then by the above inequality we get
\begin{equation*}
\begin{split}\iint_Q|T_k(u_n)|^{m'}dxdt
&\leq C2^{\sum\limits_{i=1}^{N}\frac{\bar{\rho}}{m'p_iN}+\frac{p_i-\rho_i}{p_i}\frac{\bar{\rho}}{\rho_iN}}k^{\sum\limits_{i=1}^{N}\left[\frac{\theta(p_i-1)\rho_i}{p_i-\rho_i}-m'\right]\frac{p_i-\rho_i}{p_i}\frac{\bar{\rho}}{\rho_iN}}\\
& \ \ \ \ \times\left(\iint_Q|T_k(u_n)|^{m'}dxdt\right)^{\frac{\bar{\rho}}{Nm'}+\sum\limits_{i=1}^{N}\frac{\bar{\rho}}{m'p_iN}+\frac{p_i-\rho_i}{p_i}\frac{\bar{\rho}}{\rho_iN}}\\
&\leq Ck^{\theta\bar{\rho}-\frac{\theta\bar{\rho}}{\bar{p}}-m'+\frac{m'\bar{\rho}}{\bar{p}}}
\left(\iint_Q|T_k(u_n)|^{m'}dxdt\right)^{1+\frac{\bar{\rho}}{\bar{p}m'}-\frac{\bar{\rho}}{\bar{p}}+\frac{\bar{\rho}}{Nm'}},\\
\end{split}
\end{equation*}
where we have used the fact $m'=\frac{(N+2)\bar{\rho}}{N}$.

It follows from the above inequality and (\ref{3.67}) that for $\iint_Q|T_k(u_n)|^{m'}dxdt\geq 1$ and $k\geq1$
\begin{equation}
\begin{split}\label{3.70}
&\iint_Q|T_k(u_n)|^{m'}dxdt\\
&\leq
Ck^{\left[\theta\bar{\rho}-\frac{\theta\bar{\rho}}{\bar{p}}-m'+\frac{m'\bar{\rho}}{\bar{p}}\right]\frac{\bar{p}Nm'}{(Nm'-\bar{p}-N)\bar{\rho}}}\\
&=Ck^{\frac{m[\theta N(\bar{p}-1)(m-1)-(m-1)(N+2)\bar{p}+Nm]}{(N-\bar{p}m+\bar{p})(m-1)}}.
\end{split}
\end{equation}
Recalling that $m>\frac{\bar{p}(N+2)}{\bar{p}(N+2)-N}$, by (\ref{1.1}) and (\ref{1.9}), one can deduce that
$$\frac{m[\theta N(\bar{p}-1)(m-1)-(m-1)(N+2)\bar{p}+Nm]}{(N-\bar{p}m+\bar{p})(m-1)}>0.$$
This, combined with (\ref{3.70}), gives that for any $k\geq1$
\begin{equation}
\begin{split}\label{p343}
&\iint_Q|T_k(u_n)|^{m'}dxdt\leq Ck^{\frac{m[\theta N(\bar{p}-1)(m-1)-(m-1)(N+2)\bar{p}+Nm]}{(N-\bar{p}m+\bar{p})(m-1)}}.
\end{split}
\end{equation}
If $k\leq1,$ then using (\ref{1.9}), we also have
\begin{equation}
\begin{split}\notag\iint_Q|T_k(u_n)|^{m'}dxdt\leq|Q|k^{m'} \leq|Q|k^{\frac{m[\theta N(\bar{p}-1)(m-1)-(m-1)(N+2)\bar{p}+Nm]}{(N-\bar{p}m+\bar{p})(m-1)}}.
\end{split}
\end{equation}
Combing the above two inequalities, it follows that for any $k>0$
\begin{equation}
\begin{split}\label{p344}\iint_Q|T_k(u_n)|^{m'}dxdt\leq Ck^{\frac{m[\theta N(\bar{p}-1)(m-1)-(m-1)(N+2)\bar{p}+Nm]}{(N-\bar{p}m+\bar{p})(m-1)}}.
\end{split}
\end{equation}
As a result, we obtain
\begin{equation}
\begin{split}\label{3.75}
|\left\{(x,t)\in Q: |u_n(x,t)|>k\right\}|&\leq Ck^{\frac{m[\theta N(\bar{p}-1)(m-1)-(m-1)(N+2)\bar{p}+Nm]}{(N-\bar{p}m+\bar{p})(m-1)}-m'}\\
 &= Ck^{-\frac{m[(N+2)\bar{p}-\theta N(\bar{p}-1)-N-\bar{p}]}{N-\bar{p}m+\bar{p}}}.
\end{split}
\end{equation}
Thus, we get (\ref{3.61}).

By (\ref{3.65}) and (\ref{p344}), it is easy to see that
\begin{equation}
\begin{split}\label{3.76}
\iint_Q|D_iT_k(u_n)|^{p_i}dxdt\leq C(1+k)^{\theta(p_i-1)}k^{\frac{\theta N(\bar{p}-1)(m-1)-(m-1)(N+2)\bar{p}+Nm}{N-\bar{p}m+\bar{p}}}.
\end{split}
\end{equation}
From (\ref{3.76}) and proposition 3.2, we get the desired result (\ref{3.62}). This completes the proof of case (i).

{\bf Case (ii): }$1<m\leq \frac{\bar{p}(N+2)}{\bar{p}(N+2)-N}.$

In this case, we have $m'\geq\frac{\bar{p}(N+2)}{N},$ and as a result, one obtains
\begin{equation}
\begin{split}\label{3.77}
\left(\iint_Q|T_k(u_n)|^{m'}dxdt\right)^{\frac{1}{m'}}\leq k^{1-\frac{\bar{p}(N+2)}{Nm'}}\left(\iint_Q|T_k(u_n)|^{\frac{\bar{p}(N+2)}{N}}dxdt\right)^{\frac{1}{m'}}.
\end{split}
\end{equation}
Using (\ref{3.64}), (\ref{3.65}), (\ref{3.77}) and proposition 3.1, we get
\begin{equation}
\begin{split}\label{3.80}
&\iint_Q|T_k(u_n)|^{\frac{\bar{p}(N+2)}{N}}dxdt\\
&\leq C\left(\mathop {\mbox{ess\ sup}}_{0\leq t\leq T}\int_\Omega|T_k(u_n)|^2dx\right)^{\frac{\bar{p}}{N}}
\prod_{i=1}^{N}\left(\iint_Q|D_iu_n|^{p_i}dxdt\right)^{\frac{\bar{p}}{p_iN}}\\
&\leq C\left[k^{1-\frac{\bar{p}(N+2)}{Nm'}}\left(\iint_Q|T_k(u_n)|^{\frac{\bar{p}(N+2)}{N}}dxdt\right)^{\frac{1}{m'}}\right]^{\frac{\bar{p}}{N}}\\
& \ \ \ \ \times\prod_{i=1}^{N}\left[(1+k)^{\theta(p_i-1)}k^{1-\frac{\bar{p}(N+2)}{Nm'}}\left(\iint_Q|T_k(u_n)|^{\frac{\bar{p}(N+2)}{N}}dxdt\right)^{\frac{1}{m'}}\right]^{\frac{\bar{p}}{p_iN}}\\
&=C(1+k)^{{\theta}(\bar{p}-1)}k^{\left[1-\frac{\bar{p}(N+2)}{Nm'}\right](1+\frac{\bar{p}}{N})}\left(\iint_Q|T_k(u_n)|^{\frac{\bar{p}(N+2)}{N}}dxdt\right)^{\frac{\bar{p}+N}{m'N}}.
\end{split}
\end{equation}
Since $m\leq \frac{\bar{p}(N+2)}{\bar{p}(N+2)-N},$ one can easily check that $\frac{\bar{p}+N}{m'N}\leq \frac{\bar{p}+N}{\bar{p}(N+2)}<1$ and $\frac{\bar{p}(N+2)}{Nm'}\leq 1.$  Therefore, we have
\begin{equation}
\begin{split}\label{3.81}
&\iint_Q|T_k(u_n)|^{\frac{\bar{p}(N+2)}{N}}dxdt\\
&\leq \left[C(1+k)^{\theta(\bar{p}-1)}k^{\left[1-\frac{\bar{p}(N+2)}{Nm'}\right](1+\frac{\bar{p}}{N})}\right]^{\frac{1}{1-\frac{\bar{p}+N}{m'N}}}\\
&\leq C(1+k)^{\frac{\theta mN (\bar{p}-1)}{N-\bar{p}m+\bar{p}}}k^{\frac{(N+\bar{p})[Nm-\bar{p}(m-1)(N+2)]}{N(N-\bar{p}m+\bar{p})}},
\end{split}
\end{equation}
which yields
\begin{equation}
\begin{split}\label{3.82}
&\iint_Q|T_k(u_n)|^{\frac{\bar{p}(N+2)}{N}}dxdt\\
&\leq C2^{\frac{\theta mN (\bar{p}-1)}{N-\bar{p}m+\bar{p}}}k^{{\frac{(N+\bar{p})[Nm-\bar{p}(m-1)(N+2)]+\theta mN^2 (\bar{p}-1)}{N(N-\bar{p}m+\bar{p})}}},\ \ \mbox{for}\ k\geq1 .
\end{split}
\end{equation}
By (\ref{1.1}) and (\ref{3.81}), we arrive at
\begin{equation}
\begin{split}\label{3.83}
\iint_Q|T_k(u_n)|^{\frac{\bar{p}(N+2)}{N}}dxdt&\leq |Q|k^{\frac{\bar{p}(N+2)}{N}}\notag\\
&\leq |Q|k^{\frac{(N+\bar{p})[Nm-\bar{p}(m-1)(N+2)]+\theta mN^2 (\bar{p}-1)}{N(N-\bar{p}m+\bar{p})}},\ \ \mbox{for}\ k<1.
\end{split}
\end{equation}
Combining the above inequality with (\ref{3.82}), one can deduce that
\begin{equation}
\begin{split}\label{3.84}
\iint_Q|T_k(u_n)|^{\frac{\bar{p}(N+2)}{N}}dxdt\leq Ck^{\frac{(N+\bar{p})[Nm-\bar{p}(m-1)(N+2)]+\theta mN^2 (\bar{p}-1)}{N(N-\bar{p}m+\bar{p})}}.
\end{split}
\end{equation}
This implies that (\ref{3.61}) holds true.
Finally, by (\ref{3.64}), (\ref{3.65}), (\ref{3.77}), (\ref{3.84}) and proposition 3.2, we obtain the desired result (\ref{3.62}).

{\bf Case (iii): }$m=1.$

In this case, choosing $T_k(u_n(x,t))\chi_{(0,\tau)}(t)$ as a test function for problem $(P_n)$, we have
\begin{equation}
\begin{split}\label{p352}
&\mathop {\mbox {ess\ sup}}_{0\leq t\leq T}\int_\Omega|T_k(u_n(x,\tau))|^2dx+\iint_Q\sum_{i=1}^{N}\frac{|D_iT_k(u_n)|^{p_i}}{(1+|u_n|)^{\theta(p_i-1)}}dxdt\\
&\leq k\|f\|_{L^1(Q)}
\end{split}
\end{equation}
and
\begin{equation}
\begin{split}\label{p353}
&\iint_Q|D_iT_k(u_n)|^{p_i}dxdt\leq C (1+k)^{\theta(p_i-1)}k.
\end{split}
\end{equation}
Arguing as in the proof of (\ref{3.84}), by (\ref{p352}), (\ref{p353}) and proposition 3.1, we obtain
\begin{equation}
\begin{split}\label{3.88}
\iint_Q|T_k(u_n)|^{\frac{\bar{p}(N+2)}{N}}dxdt\leq Ck^{\frac{N+\bar{p}+\theta N(\bar{p}-1)}{N}},
\end{split}
\end{equation}
which shows that (\ref{3.61}) holds true. The result (\ref{3.62}) follows from (\ref{p352})-(\ref{3.88}) and proposition 3.2 immediately.

{\bf Step 2.} In this step, we give the proof of the estimates of (\ref{3.59}) and (\ref{3.60}).

Taking $T_1(u_n){\chi(0,\tau)}(t)$ as a test function for problem $(P_n)$, from (\ref{1.2}) and H\"{o}lder's inequality, we obtain
\begin{equation}
\begin{split}\notag
\int_\Omega \int_0^{u_n(x,\tau)}T_1(s)dsdx + \alpha \int_0^\tau\int_\Omega\sum\limits_{i=1}^{N}\frac{|D_iT_1(u_n)|^{p_i}}{(1+|u_n|)^{\theta(p_i-1)}}dxdt
\leq \|f_n\|_{L^1(Q)}.
\end{split}
\end{equation}
This, together with the fact $|s|-\frac{1}{2}\leq \int_0^{s}T_1(\eta)d\eta\leq |s|,$  gives
\begin{equation}\label{3.90}
\mathop {\mbox{ess\ sup}}_{0\leq t\leq T}\int_\Omega|u_n(x,t)|dx\leq C,
\end{equation}
which is the desired result (\ref{3.59}).

From (\ref{p352}) and (\ref{p353}), it follows that (\ref{3.60}) holds true. Thus, the proof is completed.\ \ $\Box$

\section{Proofs of Theorem 2.1-2.5}

In this section, we prove the main results. Since the proofs of Theorem 2.2, Theorem 2.3 and Theorem 2.4 are very similar to that of
Theorem 2.1, we only give the proofs of Theorem 2.1 and Theorem 2.5.\vspace{2mm}

{\bf Proof of Theorem 2.1.} By Lemma 3.1 and the equation of problem $(P_n)$, we conclude that
$$\frac{\partial u_{n}}{\partial t}\ \mbox{is\ bounded\ in}\ L^{\vec{p'}}(0,T;W^{-1, \vec{p'}}(\Omega))+L^m(Q).$$
Therefore,  using a compactness result (see \cite{31}), we deduce that there exist a subsequence of $\{u_n\}$
(still denoted by $\{u_n\}$) and a measurable function $u$ such that
\begin{equation}\label{3.5}u_n\longrightarrow u \mbox{  a.e\ \ in }  Q,
 \end{equation}
 \begin{equation}\label{3.1}u_n\rightharpoonup u \mbox{ weakly\ \ in } L^{\vec{p}}(0,T;W_0^{1,
\vec{p}}(\Omega)), \end{equation}
and
\begin{equation}\label{3.3}u_n\rightharpoonup u \ \mbox{  weakly}*\ \mbox{\ in }
L^{\infty}(Q). \end{equation}
 Arguing as in the proof of Theorem 3.3 of \cite{11}, we have
 \begin{equation}\label{3.8}Du_n\longrightarrow Du \ \mbox{  a.e.\ \ in }  Q.
\end{equation}
By (\ref{1.3}), (\ref{1.4}), (\ref{3.5}), (\ref{3.8}) and using Vitali's Theorem, it follows that
\begin{equation}\label{3.9}
a_i(x,t,T_n(u_n),Du_n)\longrightarrow a_i(x,t,u,Du) \  \mbox{  strongly\ \ in }  L^{p'_i}(Q).
\end{equation}
Passing to the limit as $n\rightarrow\infty$ in the weak formulation of $(P_n)$, using (\ref{2.2}), (\ref{2.3}), (\ref{2.4}), (\ref{3.1}) and (\ref{3.9}), we infer that $u$ is a weak solution of problem $(P)$.\ \ $\Box$\vspace{2mm}

{\bf Proof of Theorem 2.5.} The proof is similar to \cite{30}, so we only sketch it here.
Taking $T_k(u_n){\chi_{(0,\tau)}}(t)$ as a test function for problem $(P_n)$, we get
\begin{equation}
\begin{split}\label{4.6}
&\int_\Omega \int_0^{u_n(x,\tau)}T_k(s)dsdx + \alpha \int_0^\tau\int_\Omega\sum\limits_{i=1}^{N}\frac{|D_iT_k(u_n)|^{p_i}}{(1+|u_n|)^{\theta(p_i-1)}}dxdt\leq k\|f\|_{L^1(Q)}.
\end{split}
\end{equation}
It can be deduced from (\ref{4.6})
and the equation of $(P_n)$ that for any $S\in W^{2,\infty}(\mathbb{R})$ such that $S'$ is compact
$$\frac{\partial S(u_{n})}{\partial t}\ \mbox{is\ bounded\ in}\ L^{\vec{p'}}(0,T;W^{-1, \vec{p'}}(\Omega))+L^1(Q).$$
 Arguing again as in \cite{Blan98}, the above two estimates imply that there exist a subsequence of $\{u_n\}$
(still denoted by $\{u_n\}$) and a measurable function $u$ such that
\begin{equation}\label{4.7}u_n\longrightarrow u \mbox{  a.e\ \ in }  Q,
 \end{equation}
 and
\begin{equation}\label{4.31}
T_k(u_n)\rightharpoonup T_k(u)  \ \mbox { weakly in } L^{\vec{p}}(0,T;W_0^{1,\vec{p}}(\Omega)).
\end{equation}
Then, proceeding as in \cite{30,22}, we get
 \begin{equation}\label{4.8} Du_n\longrightarrow  Du \ \mbox{  a.e.\ \ in }  Q.
\end{equation}
Arguing as in Theorem 3.2 of \cite{30}, we obtain that $u\in L^\infty(0,T;L^1(\Omega))$ is a renormalized solution of $(P)$, i.e. $u$ satisfies
\begin{equation}T_k(u)\in L^{\vec{p}}(0,T;W_0^{1,\vec{p}}(\Omega)),\ \ S_k(u)\in C([0,T];L^1(\Omega)),
\end{equation}
\begin{equation}
\begin{split}\label{4.9}
&\underset{n\rightarrow\infty}{\mbox{lim}}\iint_{\{n\leq |u|\leq n+1\}}\sum\limits_{i=1}^{N}a_i(x,t,u,Du)D_iudxdt=0
\end{split}
\end{equation}
and
\begin{equation}
\begin{split}\label{4.10}
&(S(u))_t-\mbox{div}\left(S'(u)a(x,t,u,Du)\right)\\
&+\sum\limits_{i=1}^{N}a_i(x,t,u,Du)D_iuS''(u)=fS'(u),
\end{split}
\end{equation}
where $a(s,t,u,Du)=(a_1(x,t,u,Du),a_2(x,t,u,Du),\ldots,a_N(x,t,u,Du))$.

Let us choose $S(s)=H_n(s)=\int_0^s 1-|T_1(\tau-T_n(\tau))|d\tau$, and take $T_k(H_n(u)-\phi)$ as test function in (\ref{4.10}).  Passing to
the limit as $n\rightarrow\infty$, one can deduce that $u$ is a entropy solution of problem $(P)$.

Furthermore, by (\ref{3.61}), (\ref{3.62}) and Fatou's lemma, we
get $u\in {\cal M}^{\tilde{r}}(Q)$ and $|D_iu|\in {\cal M}^{\beta_i}(Q).$
 Hence the proof of Theorem 2.5 is completed.\ \ $\Box$\vspace{2mm}

\noindent\textbf{Remark 4.1}  Without loss of generality, we may assume that $p_1\geq p_2\geq p_3\geq ...\geq p_N$.
If $f\in L^m(\xo)$ with $m_2:=\frac{\bar{p}[(2\theta-N-2)p_1-\theta p_2(p_1+1)]}{\bar{p}[(2\theta-N-2)p_1-\theta p_2(p_1+1)]+Np_1[1+\theta(p_2-1)]}\leq m<m_1,$ we have $1+\theta(p_2-1)\leq l<1+\theta(p_1-1).$ Arguing as in the proof of Lemma 3.4, we conclude that problem $(P_n)$ admits at least a weak solution $u_n\in L^{\vec{q}}(0,T;W_0^{1,\vec{q}}(\Omega))\cap
L^{\tilde{r}}(Q)$ with $\vec{q}=(q_1,q_2,...,q_N)$ and $\tilde{r}=\frac{N+l+1}{N}\bar{q}=lm'$ satisfying
$$\|u_n\|_{L^{\vec{q}}(0,T;W_0^{1,\vec{q}}(\Omega))\cap
L^{\tilde{r}}(Q)}\leq C,$$
where
 $$q_i=\left\{\begin{array}{ll}
\frac{(N+l+1)\bar{p}p_1-\bar{p}p_1\theta-\bar{p}(1-l-\theta)}{\bar{p}p_1(N+l+1)-\bar{p}(p_1\theta+1-l-\theta)+Np_1[\theta(p_1-1)+1-l]} & i=1,\\
 p_i &
i\geq 2,\\
 \end{array}
 \right.$$
 and
$$l=\frac{\bar{p}\left[\theta+(1-\theta)\frac{1}{p_1}-N-1\right]}{\bar{p}(1+\frac{1}{p_1})-m'N}.$$
Similar to the proof of Theorem 2.1, we can infer that there exists at least a solution $u\in L^{\vec{q}}(0,T;W_0^{1,\vec{q}}(\Omega))\cap
L^{\tilde{r}}(Q)$ to problem $(P)$, where $\vec{q}$ and $\tilde{r}$ are defined as above.

If $f\in L^m(\xo)$ with $m_3:=\frac{\bar{p}[p_1p_2(2\theta-N-2)-\theta(p_1+p_2)p_3]}{\bar{p}[p_1p_2(2\theta-N-2)-\theta(p_1+p_2)p_3]+Np_1p_2[1+\theta(p_3-1)]}\leq m<m_2,$
we have $1+\theta(p_3-1)\leq l< 1+\theta(p_2-1).$ Arguing as before, we infer that problem $(P)$ admits at least a weak solution $u\in L^{\vec{q}}(0,T;W_0^{1,\vec{q}}(\Omega))\cap
L^{\tilde{r}}(Q)$ with $\vec{q}=(q_1,q_2,...,q_N)$ and $\tilde{r}=\frac{N+l+1}{N}\bar{q}=lm'$, where
 $$q_i=\left\{\begin{array}{ll}
\frac{\bar{p}p_ip_1p_2(N+l+1)-\bar{p}p_i[2\theta p_1p_2+(1-l-\theta)(p_1+p_2)]}{Np_1p_2[\theta(p_i-1)+1-l]+p_1p_2\bar{p}(N+l+1)-\bar{p}[2\theta p_1p_2+(1-l-\theta)(p_1+p_2)]} & i=1,2,\\
 p_i &
i\geq 3,\\
 \end{array}
 \right.$$
 and $l=\frac{\bar{p}\left[2\theta+(1-\theta)(\frac{1}{p_1}+\frac{1}{p_2})-N-1\right]}{\bar{p}(1+\frac{1}{p_1})+\frac{1}{p_2}-m'N}.$ Thus, using exactly the same procedure presented as above, we shall divide this interval $[m_N, m_1)$ into $N-1$ subintervals and obtain the corresponding results for each subintervals. For the last interval, we have the following result.

If $f\in L^m(\xo)$ with $m_N\leq m
<m_{N-1}$,
we have $1+\theta\left(\min\limits_{1\leq i\leq N}p_i-1\right)\leq l< 1+\theta(p_{N-1}-1),$ where $m_{N-1}=\frac{\bar{p}[(N\theta-N-2+\theta)\bar{p}p_N-\theta p_{N-1}(Np_N-\bar{p}+\bar{p}p_N)]}{\bar{p}[(N\theta-N-2+\theta)\bar{p}p_N-\theta p_{N-1}(Np_N-\bar{p}+\bar{p}p_N)]+N\bar{p}p_N[1+\theta (p_{N-1}-1)]}.$ In this case, we may conclude that problem $(P)$ admits at least a weak solution $u\in L^{\vec{q}}(0,T;W_0^{1,\vec{q}}(\Omega))\cap
L^{\tilde{r}}(Q)$ with $\vec{q}=(q_1,q_2,...,q_N)$ and $\tilde{r}=\frac{N+l+1}{N}\bar{q}=lm'$, where
 $$q_i=\left\{\begin{array}{ll}
 \frac{(N+l+1)\bar{p}p_i-\bar{p}p_i\left[(N-1)\theta+(1-l-\theta)(\frac{N}{\bar{p}}-\frac{1}{p_N})\right]}{N[\theta(p_i-1)+1-l]+(N+l+1)\bar{p}-\bar{p}\left[(N-1)\theta+(1-l-\theta)\left(\frac{N}{\bar{p}}-\frac{1}{p_N}\right)\right]}& i=1,2,...,N-1,\\
 p_i &
i=N,\\
 \end{array}
 \right.$$
and $l=\frac{\bar{p}\left[N\theta+(1-\theta)\left(\frac{N}{\bar{p}}-\frac{1}{p_N}\right)-N-1\right]}
{\bar{p}(1+\frac{N}{\bar{p}}-\frac{1}{p_N})-m'N}.$\\

\textbf{Acknowledgment}

This research is supported by National Natural Science Foundation of China(No.11801259, No.11461048), Natural Science Foundation of Jiangxi Province (No.20202BABL201009, No. 20202ACBL211002), Academic and Technical Leaders Training Plan of Jiangxi Province (No. 20212BCJ23027) and the foundation of Education Department of Jiangxi
Province(no.GJJ190514).

 \small{
}

  \end{document}